\documentclass[11pt]{article}
\usepackage{epsfig,graphics,graphicx}
\usepackage{amssymb}
\usepackage{amsmath}
\usepackage{extarrows}
\usepackage[utf8]{inputenc}
\usepackage[english]{babel}
\usepackage{geometry}
\usepackage{amsthm}
\usepackage{scalerel}[2014/03/10]
\usepackage[usestackEOL]{stackengine}
\usepackage{esint}
\usepackage{xcolor}

\newtheorem{theo}{Theorem} 

\newtheorem{lem}[theo]{Lemma}
\newtheorem{prop}[theo]{Proposition}
\newtheorem{coro}[theo]{Corollary}

\def\R{\mathbb R}

\def\C{\mathbb C}

\def\bC{\mathbf C}

\def\cK{\mathcal K}

\def\bM{\mathbf M}

\def\Re{\text{Re}}
\def\Im{\text{Im}}

\def\cB{\mathcal B}

\def\diam{\text{diam}}
\def\bv{\mathbf{v}}

\def\div{\operatorname{div}}
\def\Div{\operatorname{div}}
\def\supp{\operatorname{supp}}

\def\Id{\operatorname{\textbf{Id}}}

\def\curl{\operatorname{curl}}

\def\v{\operatorname{\textbf{v}}}

\newcommand{\norm}[2][]

\title{\sc{Well-posedness of compressible $2D$ \\nonlinear transport equations\\ and quasiconformal maps}}
\author{
\sc{J.C. Cantero}\thanks{%
    Dpt. Matematics, Universitat Aut\`onoma de Barcelona, 
    Email: \textit{juancarlos.cantero@uab.cat}, }
\sc{A. Clop}\thanks{%
    Dpt. Mathematics and Computer Science, Universitat de Barcelona, 
    Email: \textit{albert.clop@ub.edu}, }
\sc{D. Faraco}\thanks{%
    Department of Mathematics, Universidad Aut\`onoma de Madrid, 
    Email: \textit{daniel.faraco@uam.es}, }
\sc{D. Horas}\thanks{%
    Dpt. Mathematics and Computer Science, Universitat de Barcelona, 
    Email: \textit{david.horas@ub.edu}, }, \\
\sc{J. Mateu}\thanks{%
    Dpt. Mathematics, Universitat Aut\`onoma de Barcelona, 
    Email: \textit{joan.mateu@uab.cat}, }
\sc{J. Orobitg}\thanks{%
    Dpt. Mathematics, Universitat Aut\`onoma de Barcelona, 
    Email: \textit{joan.orobitg@uab.cat}, }
\sc{B. Sengupta} \thanks{%
    Centre for Applicable Mathematics, Tata Institute of Fundamental Research,
Email: \textit{sengupta25@tifrbng.res.in}, } 
}
\date{}

\author{%
  \sc{Juan Carlos Cantero}
  \thanks{
  \texttt{juancarlos.cantero@uab.cat}, \texttt{joan.mateu@uab.cat}, \texttt{joan.orobitg@uab.cat};
    Departament de Matemàtiques, Universitat Autònoma de Barcelona, 08193-Cerdanyola del Vall\`es, Catalonia.
  }
  \and
  \sc{Albert Clop}
  \thanks{
  \texttt{albert.clop@ub.edu}, \texttt{david.horas@ub.edu};
    Department of Mathematics and Computer Science, Universitat de Barcelona, 08007-Barcelona, Catalonia.
  } 
\and
  \sc{Daniel Faraco}
  \thanks{
  \texttt{daniel.faraco@uam.es};
    Department of Mathematics, Universidad Aut\'onoma de Madrid, 28049-Madrid, Spain.  
  }  
  \and
  \sc{David Horas} \footnotemark[2]
  \and
  \sc{Joan Mateu} \footnotemark[1]
  \and
  \sc{Joan Orobitg} \footnotemark[1]
  \and
  \sc{Banhirup Sengupta} \thanks{
  \texttt{sengupta25@tifrbng.res.in}; Centre for Applicable Mathematics, Tata Institute of Fundamental Research; Bangalore 560065, India.
  }
}

\begin{document}

\maketitle

\abstract{We prove uniqueness of solutions to a class of nonlinear transport equations in the plane, for which the velocity field is obtained as the convolution of the classical Cauchy kernel with the unknown. Even though the initial datum is bounded and compactly supported, the velocity field may have unbounded divergence. The proof relies on the theory of quasiconformal mappings.}

\section{Introduction}

\noindent
We study the solvability of the following initial value problem
\begin{equation}\label{cauchysystem}
\begin{cases}
\frac{d}{dt}\omega+\v\cdot \nabla \omega = 0, \\
\v(t,\cdot)=K_\theta\ast \omega(t,\cdot),\\
\omega(0,\cdot)=\omega_0.
\end{cases}\end{equation}
where $K_\theta$ is the following convolution kernel 
$$\aligned
K_\theta(x,y) =\frac1{2\pi}\,\frac{(x\,\cos\theta+y\,\sin\theta, x\,\sin\theta-y\,\cos\theta) }{x^2+y^2}=\frac{e^{i\theta}}{2\pi z},
\endaligned$$
and $\theta\in [0, 2\pi]$ is fixed. The easiest situation $\theta=0$ corresponds to the classical Cauchy kernel,
$$K_0(z)=\frac1{2\pi z}=\frac1{2\pi}\,\frac{(x,-y)}{x^2+y^2}$$ 
with the identification $\R^2=\C$. The datum $\omega_0\in L^\infty$ is a given compactly supported and real-valued function, and $\v$ is the velocity field. A weak solution with initial datum $\omega_0$ is a map $\omega:[0,T]\times\R^2\to\R$  with $\omega\in L^\infty\bigl([0,T];L^1 \cap L^\infty \bigr)$ and such that
\begin{equation}\label{eq:weak}
\int_0^T\!\!\int_{\R^2}
\omega\bigl(\varphi_t+\mathbf{v}\cdot\nabla\varphi+(\Div \mathbf{v})\varphi\bigr) \,dA\,dt
+\int_{\R^2}\omega_0\varphi(0) \,dA=0.
\end{equation}
for every \(\varphi\in C_c^\infty([0,T)\times\R^2)\). The assumption $\omega\in L^\infty([0,T];L^1 \cap L^\infty )$ suffices to guarantee that the integrals in \eqref{eq:weak} are absolutely convergent. The model \eqref{cauchysystem} is a natural counterpart to the \emph{planar Euler system of equations}, which in its vorticity form takes the form \eqref{cauchysystem} with the kernel
$$\aligned
K(z)=\frac{i}{2\pi\,\bar{z}}=\frac1{2\pi}\,\frac{(-y,x)}{ x^2+y^2}.
\endaligned$$
In the same way Euler's system describes the time evolution of a perfect, incompressible, inviscid fluid, with vorticity $\omega$, the model with kernel $K_\theta$ can be interpreted, at least formally, as the transport of an internal activity field by the flow it generates in an idealized two-dimensional compressible viscous medium with purely deviatoric stresses and fixed preferred axes oriented at angles $\theta/2$ and $\theta/2+\pi/2$ to the horizontal. 

\vspace{.5cm}

Our main result is that, if $\theta\in[0,2\pi]$ is fixed,  the problem \eqref{cauchysystem} is well posed for any bounded and compactly supported initial datum $\omega_0$.

\vspace{.5cm}

\begin{theo}\label{main1}
The initial value problem \eqref{cauchysystem} is well posed for each initial datum $\omega_{0}\in L^\infty$ with compact support. That is, for every $T>0$ there is exactly one weak solution $\omega\in L^\infty([0,T];L^1\cap L^\infty)$. 
 \end{theo}

\vspace{.5cm}

Existence of weak solutions $\omega\in L^\infty([0,T];L^1\cap L^\infty)$ was proven in \cite{CS} for compactly supported $\omega_0\in L^\infty$. However, uniqueness remained an open question which is solved here.  As explained in \cite{CS}, the choice of the kernel $K_\theta$ ensures that
\begin{equation}\label{dbarv}
\overline\partial\v =\frac{e^{i\theta}\,\omega}{2}
\end{equation}
where $ \partial=\frac{\partial_x-i\partial_y}{2}$ and $\overline\partial=\frac{\partial_x+i\partial_y}{2}$ are the classical Wirtinger derivatives. Since $\omega\in L^\infty([0,T];L^1\cap L^\infty)$, then endpoint Calder\'on-Zygmund estimates give that $D\v\in BMO$. More remarkably, $\div(\v)$ may not be identically $0$, thus making $\v$ compressible. In fact, even assuming that $\omega_0\in L^\infty$ is compactly supported, one could have unbounded $\div(\v)$. Thus, even though  $\v$ admits a well-defined flow $X:[0,T]\times\R^2\to\R^2$ of H\"older continuous homeomorphisms $X(t,\cdot)$, they need not be measure preserving, thus putting at risk the control of both $L^1$ and $L^\infty$ norms. In this context, it turns out that certain linear transport equations \cite{CJMO, CJMO2, CJMO3} still exhibit nice well-posedness results even having unbounded divergence drifts. The main explanation is that, in that situation, one still has a well-defined flow (in the DiPerna-Lions sense, see \cite{DiPL}) that preserves Lebesgue null sets. We rely precisely on this fact to conclude that all weak solutions of \eqref{cauchysystem} are indeed Lagrangian.  

\vspace{.5cm}

As a second consequence of \eqref{dbarv}, one has $\overline\partial\v \in L^\infty([0,T];L^\infty)$ for any weak solution $\omega\in L^\infty([0,T];L^1\cap L^\infty)$. Under this assumption, Reimann's Theorem \cite{R} guarantees the quasiconformality of all $t$-advance maps $X(t,\cdot):\R^2\to\R^2$ of the flow of $\v$. In the same way quasiconformality was essential for the existence argument in \cite{CS}, the proof of uniqueness we are presenting here also relies on it. For instance, it is precisely the quasiconformality of $X(t,\cdot)$ that guarantees a good control of the integral means of Jacobian determinants (see Lemma \ref{linearinq} below).

\vspace{.5cm}

The classical scheme of 2D Euler's uniqueness does not work in our setting. The reason is that in our setting the difference $\v_1-\v_2$ of two drifts (with the same initial datum $\omega_0$) fails to belong to $L^2$, and this is due precisely  to the lack of incompressibility. Instead, a more suitable energy to be considered here is given in terms of the anticonformal derivative $\overline\partial X(t,\cdot)$ of the flow. Moreover, a slight modification of the uniqueness proof allows us to obtain the following stability result.

\vspace{.5cm}
\begin{theo}\label{main2}
Suppose that $\omega_{01},\omega_{02}\in L^\infty$ have compact support inside a disk $D\subset\C$. Let $\omega_1,\omega_2\in L^\infty([0,T];L^1\cap L^\infty)$ be the  weak solutions to \eqref{cauchysystem}, and let $X_1(t), X_2(t)$ be their corresponding flows. Let $M=\max\{\|\omega_{01}\|_{L^\infty},\|\omega_{02}\|_{L^\infty}\}$. Given $k_0\in(0,1)$, choose any $\gamma_0>2\frac{1+k_0}{1-k_0}$ and let $\tau=\frac1{M}\,\log(\frac{1+k_0}{1-k_0})$. Then there exists a constant $C=C(k_0,\gamma_0,M,D)$ such that
$$\|\overline\partial X_1 (t)-\overline\partial X_2(t)\|_{L^2}\leq C\,\|\omega_{01}-\omega_{02}\|_{L^1}^{\frac{\theta(t)}{\gamma_0}}\,\|\omega_{01}-\omega_{02}\|_{L^\infty}^{\theta(t)(1-\frac1{\gamma_0})} $$
for any $t\in[0,\min\{T,\tau\}]$, where $\theta(t)=e^{-Ct}$.
 \end{theo}
 \noindent

The paper is structured as follows. In Section \ref{preliminaries} we show that every weak solution of \eqref{cauchysystem} is indeed a Lagrangian solution, and use this fact to obtain uniform control over the corresponding Lagrangian flow. Section \ref{qcsection} is devoted to exploiting the quasiconformal nature of the flow, which is used later in \ref{energy} to obtain an energy estimate to control the difference of anticonformal derivatives of the flows associated with two weak solutions. Section \ref{proofofuniqueness} is devoted to proving Theorem \ref{main1} from the energy estimate. Last, Section \ref{stability} presents a stability result whose proof can be obtained with minor modifications in the proof of uniqueness. Unless otherwise stated, all $L^p$ spaces and norms are understood on $\C=\R^2$ with respect to the planar Lebesgue measure $dA$.
\vspace{1cm}

{\bf{Acknowledgements}}. ChatGPT Astra was used in this work. In particular, it suggested deriving an energy inequality in terms of $\overline \partial X$. The idea was merged with the knowledge of the authors to produce this manuscript. The paper was humanly written. 

AC and DH were partially supported by projects 2025-SGR-0075 (Catalonia) and PID2024-156429NB-I00 (Spain).
JCC, JM and JO were partially supported by project PID2024-155320NB-I00 (Spain). DF was partially supported by project PID2024-158664NB-C22 (Spain) and by the Severo Ochoa Programme for Centers of Excellence Grant CEX2019-000904-S funded by MCIN/AEI/ 10.13039/501100011033.
BS  would like to thank the Department of Atomic Energy, Government of India, for financial assistance and the Tata Institute of Fundamental Research, Centre for Applicable Mathematics (TIFR-CAM) for providing a stimulating scientific environment and
resources. All authors acknowledge support from Grants RED2022-134784-T and RED2018-102650-T funded by MCIN/AEI/10.13039/501100011033.

\section{Uniform estimates}\label{preliminaries}

It is well known that the distributional derivatives of the integral kernel $K_\theta(z)=\frac{e^{i\theta}}{2\pi z}$ are given by a Dirac mass and a Cauchy principal value,
\begin{equation}\label{distderiv}
    \overline\partial K_\theta=\frac{e^{i\theta}}{2}\,\delta_0,\hspace{2cm}\partial K_\theta= -\frac{e^{i\theta}}{2}\,PV\frac1{\pi z^2}.
\end{equation} 
Above, the principal value is the limit $PV\frac1{\pi z^2}=\lim_{\epsilon\to 0}\frac{\chi_{|z|>\epsilon}}{{\pi z^2}}$ in the sense of distributions. One usually calls $\cB(f)=-PV\frac{1}{\pi z^2}\ast f$ the Beurling transform of $f$, whenever the convolution makes sense. Recall that $\cB$ is an even Calder\'on-Zygmund singular integral operator, characterized by the fact that $\cB(\overline\partial\varphi)=\partial\varphi$ for any $\varphi\in C^\infty_c(\C)$. Thus $\cB$ can be represented by the unimodular Fourier multiplier $\frac{\xi}{\bar \xi}$. In particular, $\cB$ is an isometry in $L^2(\C)$, and moreover it extends boundedly to $L^p$, $1<p<\infty$, with norm recently found to be $\|\cB\|_p=p-1$ for $p>2$, see \cite{AFBeurling}.

\begin{lem}\label{LinftyK}
Let $K:\C\setminus\{0\}\to\C$ be such that $|K(z)|\leq \frac{C}{|z|}$ for some constant $C>0$.  Let $f:\C\to\C$ be given, and assume that $f\in L^1\cap L^\infty$. 
\begin{itemize}
\item[(a)] Then $K\ast f\in L^\infty$ and 
\begin{equation}\label{without support}
\|K\ast f\|_{L^\infty}\leq C\,\|f\|_1^{\frac12}\,\|f\|_{L^\infty}^{\frac12},
\end{equation} 
\item[(b)] If $f$ is compactly supported, then 
$$\limsup_{|x|\to\infty} |x|\,|K\ast f(x)|\leq C<+\infty$$
with $C$ depending only on $\|f\|_{L^\infty}$ and $|\supp(f)|$. 
\end{itemize}
\end{lem}
\begin{proof}
Let us consider a real number $R>0$. For each such $R$ and given any two arbitrary points $x$ and $y$ in the plane, we can always divide the plane into two regions, $|x-y|\leq R$ and $|x-y|>R$. Therefore,
$$\aligned
|K\ast f(x)|
&\leq\int_{\C}|K(x-y)\,f(y)|\,dA(y)\\
&=\int_{|x-y|\leq R}|K(x-y)\,f(y)|\,dA(y)+\int_{|x-y|>R}|K(x-y)\,f(y)|\,dA(y)\\
&\leq\|f\|_{L^\infty}\,\int_{|x-y|\leq R}\frac{C}{|x-y|}\,dA(y)+\frac{C}{R}\int_{|x-y|>R}|f(y)|dy\\
&\leq CR\,\|f\|_{\infty}+\frac{C}{R}\|f\|_1
\endaligned$$ 
If we minimize the term on the right-hand side of the inequality as a function of $R$, the best possible value attainable is $R=\|f\|_{L^1}^\frac12\,\|f\|_{L^\infty}^{-\frac12}.$ This gives the bound \eqref{without support}.  
Concerning (b), let us assume that $diam(\supp(f))=2R$. It is not restrictive to assume $0\in\supp(f)$, so that $\supp(f)\subset \overline{D(0,R)}$. Then, at points $x$ such that $|x|>2R$ one has
$$\aligned
|K\ast f(x)|
&\leq C\int_{|y-x|\leq R}\frac{|f(x-y)|}{|y|}\,dA(y)\\ 
&\leq C\|f\|_{L^\infty}\,\int_{|y-x|\leq R}\frac{1}{|y|}\,dA(y)\\ 
&\leq C\|f\|_{L^\infty}\,\frac{R^2}{|x|-R}\leq C\,\frac{\|f\|_{L^\infty}\,\diam(\supp(f))^2}{|x|}, 
\endaligned$$
as claimed.
\end{proof}

We now apply the above estimates to the family of kernels $K_\theta(z)=\frac{e^{i\theta}}{2\pi z}$.  By definition, any weak solution  $\omega\in L^\infty([0,T];L^1\cap L^\infty)$ of \eqref{cauchysystem} induces a drift $\v$ through the rule $\v(t)=K_\theta\ast \omega(t)$. The size of $K_\theta$ allows us to apply Lemma \ref{LinftyK}, which immediately provides us with uniform bounds for $\v$,
$$\sup_{t\in[0,T]}\|\v(t) \|_{  L^\infty }\leq C  \sup_{t\in[0,T]}\|\omega(t)\|_{L^1}^\frac12\,\|\omega(t)\|_{L^\infty}^\frac12 
.$$
Also, \eqref{distderiv} gives
$$\overline\partial \v(t)=\frac{e^{i\theta}}{2}\,\omega(t)\qquad\partial\v(t)=\frac{e^{i\theta}}{2}\,\cB(\omega(t)).$$ 
which combined with the $L^p$ estimates for $\cB$ give us that 
$$\aligned
 \sup_{t\in[0,T]}\|D\v(t)\|_{L^p}&= \sup_{t\in[0,T]}\frac12\,\||\omega(t)|+|\cB\omega(t)|\|_{L^p}
 \leq  \sup_{t\in[0,T]}\frac{p}2  \|\omega(t)\|_{L^p}
 \leq  \sup_{t\in[0,T]}\frac{p}2\, \|\omega(t)\|_{L^1}^\frac1p\,\|\omega(t)\|_{L^\infty}^{1-\frac1p}
\endaligned$$
whenever $1<p<\infty$. Moreover, the endpoint estimate $\cB:L^\infty(\C)\to BMO(\C)$ also gives that
\begin{equation}\label{eq:BMO}
 \sup_{t\in[0,T]}\|D\v(t)\|_{BMO}\leq C\, \sup_{t\in[0,T]}\|\omega(t)\|_{L^\infty}.
 \end{equation}
From \eqref{eq:BMO}, exponential bounds for $\div(\v)$ follow immediately, due to the John-Nirenberg inequality. This fact puts $\v$ into the theory initiated in \cite{DJ} and further developed later in \cite{CJMO, CJMO2} for vector fields with exponentially integrable divergence (not directly covered by the classical DiPerna-Lions theory from \cite{DiPL}). Roughly speaking, initial value problems for the linear transport equation $u_t + \v\cdot\nabla u=0$ are well posed, and moreover their solutions are indeed Lagrangian, that is, they can be recovered from the Lagrangian flow of $\v$ (see \cite{CJMO2}). 

\vspace{.5cm}

\begin{theo}\label{weakilagrangian}
Let $\omega\in L^\infty([0,T];L^1\cap L^\infty)$ be a weak solution of \eqref{cauchysystem} with $\omega_0\in L^\infty$. Then,  
$$\omega(t, X(t,x))=\omega_0(x)$$
where $X:[0,T]\times\C\to\C$ is the DiPerna-Lions flow of ${\mathrm{\mathbf{v}}}$. Moreover, $X(t,\cdot)$ preserves Lebesgue null sets.
\end{theo}
\begin{proof}
The weak solution $\omega$ produces a velocity field $\v=K_\theta\ast\omega$. It follows from \eqref{eq:BMO} that $\v(t,\cdot)$  has 
BMO distributional derivatives and therefore 
is a Lip-log continuous map. Then, we know from \cite[Theorem 1]{CJMO} that the linear transport equation
\begin{equation}\label{transport}
\begin{cases}
    \rho_t+\v\cdot\nabla\rho=0\\\rho(0,\cdot)=\omega_0
\end{cases}
\end{equation}
admits a unique weak solution $\rho$. This solution can only be $\rho=\omega$. On the other hand, the Lip-log nature of $\v$ ensures that the drift $\v$ admits a unique well-defined flow $X:[0,T]\times \C\to\C$, and moreover for every $t\in[0,T]$ the map $X(t,\cdot):\C\to\C$ has Lusin's N property, that is, if $E\subset\C$ has zero Lebesgue measure then also its image set $X(t,E)$  does \cite[Main Theorem]{CJMO2}. Hence, the composition $\omega_0(X^{-1}(t,\cdot))$ is a well-defined, measurable and essentially bounded function. Then it just remains to see that $\omega_0(X^{-1}(t,\cdot))$ is indeed a weak solution of \eqref{transport}: to do this, just test $\omega_0(X^{-1}(t,\cdot))$ on \eqref{eq:weak} against a smooth, compactly supported function $\varphi$ on $[0,T)\times\C$. A change of variables and an integration by parts both in time and space separately (as explained in \cite{CJMO3}) will immediately show that \eqref{eq:weak} is satisfied. Thus, using uniqueness again,  necessarily $\omega_0(X^{-1}(t,\cdot))=\omega$ as claimed.
\end{proof}

In particular, every weak solution $\omega\in L^\infty([0,T];L^1\cap L^\infty)$ is generated by the flow of its corresponding drift $\v=K_\theta\ast\omega$,
$$\begin{cases}
    \frac{d}{dt}\,X(t,z)=\v(t,X(t,z)),\\X(0,z)=z.
\end{cases}$$
The existence and main properties of this flow $X$ were proven in \cite{CJMO2}. This, together with the fact that $\v\in L^\infty([0,T];L^\infty)$, ensures uniform bounds for $X(t,\cdot)$. 

\quad

\begin{coro}\label{boundeddispl}
Let $\omega\in L^\infty([0,T];L^1\cap L^\infty)$ be a weak solution of \eqref{cauchysystem} with compactly supported datum $\omega_0\in L^\infty$. Let ${\mathrm{\mathbf{v}}}=K_\theta\ast \omega$, and let $X$ be the flow of ${\mathrm{\mathbf{v}}}$.
\begin{itemize}
\item[(a)] If $|x|\leq r$ then $|X(t,x)|\leq r+T\|{\mathrm{\mathbf{v}}}\|_{L^\infty([0,T];L^\infty)}.$
\item[(b)] If $\supp(\omega_0)\subset \overline{D(0,r)}$ and $0\leq t\leq T$ then $\supp(\omega(t))\subset\overline{D(0, r')}$ for some $r'\leq r+T\|{\mathrm{\mathbf{v}}}\|_{L^\infty([0,T];L^\infty)}$.
\item[(c)] There is $C>0, R>0$ such that $\sup_{0\leq t\leq T}|X(t,x)-x|\leq \frac{C}{|x|}$ for each $|x|>R$.
\end{itemize}
\end{coro}
\begin{proof}
From Theorem \ref{weakilagrangian}, we know that $\omega$ has the representation $\omega(t,x)=\omega_0(X^{-1}(t,x))$, where $X$ is the flow of $\v=K_\theta\ast\omega$. Thus $X$ admits an integral representation
$$X(t,x)=x+\int_0^t \v(s,X(s,x))\,ds.$$
From  Theorem \ref{weakilagrangian}, $X(s,\cdot)$ preserves Lebesgue null sets, and therefore $|\v(s,X(s,x))|\le \|\v(s,\cdot)\|_{L^\infty }$. Thus
$$\sup_{x}|X(t,x)-x|\leq T\|\v\|_{L^\infty([0,T];L^\infty)},$$
and therefore if $|x|\leq r$ one immediately gets (a). Next, (b) is immediate from (a) since 
$\supp(\omega(t))\subset X(t,D(0,r))$. In order to show (c), we take from (b) a real number $R_0>0$ such that $\supp(\omega(t))\subset \overline{D(0,R_0)}$ for each $t\in[0,T]$. For instance, $R_0=r+T\|{\mathrm{\mathbf{v}}}\|_{L^\infty([0,T];L^\infty)}$ would do the job if $\supp(\omega_0)\subset D(0,r)$. Next, we observe from the integral equation that 
$$|X(t,x)|\geq |x|-T\|\v\|_{L^\infty([0,T];L^\infty)}\hspace{2cm}t\in[0,T],\,\,x\in\C.$$
Thus, at points $x$ with $|x|\geq 2R_0 +T\|\v\|_{L^\infty([0,T];L^\infty)}$ we have $|X(s,x)|\geq 2R_0 $ for every $s\in[0,T]$. Therefore from Lemma \ref{LinftyK} (b) we get at those points $x$ that
$$|\v(s, X(s,x))|\leq \frac{C\,\|\omega(s)\|_{L^\infty}\,R_0^2}{|X(s,x)|},\qquad 0\leq s\leq T$$
whence
$$|X(t,x)-x|
\leq \frac{C\,R_0^2}{|x|-T\|\v\|_{L^\infty([0,T];L^\infty)}}\,\int_0^T \|\omega(s)\|_{L^\infty}\,ds.$$
Thus (c) follows after choosing $R=2R_0+T\|\v\|_{L^\infty([0,T];L^\infty)}$, for some constant $C$ depending only on $\|\omega\|_{L^\infty([0,T];L^1\cap L^\infty)}$ and $T$.
\end{proof}

\section{Quasiconformal flows}\label{qcsection}

Given $K\geq 1$, let us recall here that a homeomorphism $f:\C\to\C$ is $K$-quasiconformal if $f\in W^{1,2}_{loc}(\C;\C)$ and moreover $|Df(z)|^2\leq K\,\det Df(z)$ at almost every $z\in\C$. We refer the reader to the monograph \cite{AIM} to understand the role of quasiconformal maps in several independent problems. In our context, quasiconformal maps arise in connection with the following celebrated result, due to H.M. Reimann, see \cite{R}.  We only state it on the plane, although it holds also in higher dimensions.

 \vspace{.5cm}
\begin{theo}[{\bf{Reimann}}]\label{Reimannth}
Let $\bv:[0,T]\times\C\to\C$ be a continuous vector field, such that for each $t$ one has
$$\limsup_{|z|\to+\infty}\frac{|\bv(t,z)|}{|z|\,\log|z|}<+\infty.$$
Suppose that the distributional derivatives $\partial \bv(t,\cdot)$ and $\overline\partial \bv(t,\cdot)$ are locally integrable functions of $z\in\C$, and moreover suppose that 
$$\sup_{t\in[0,T]}\|\overline\partial\bv (t,\cdot)\|_{L^\infty} \leq C_0<\infty.$$
Then, $\bv$ admits a unique flow $X(t,z)$ of $\cK_t$-quasiconformal maps $X(t,\cdot):\C\to\C$, and
$$\cK_t\leq \exp \left(2\int_0^t \|\overline\partial \bv(s,\cdot)\|_{L^\infty} \,ds\right)$$
whenever $0\leq t\leq T$.
\end{theo}
 \vspace{.5cm}

Concerning our initial value problem \eqref{cauchysystem}, Reimann's Theorem has the following consequence.
 \vspace{.5cm}

\begin{coro}\label{qcflows}
Let $\omega\in L^1([0,T];L^1\cap L^\infty)$ be a weak solution of \eqref{cauchysystem} with compactly supported datum $\omega_0$. Then ${\mathrm{\mathbf{v}}}=K_\theta\ast\omega$ has a well-defined flow $X:[0,T]\times\C\to\C$, and for every $t\in[0,T]$ the map $X(t,\cdot):\C\to\C$ is $\cK_t$-quasiconformal, with $1\leq \cK_t\leq e^{|t|\|\omega_0\|_{L^\infty}}$.

\end{coro}
\begin{proof}
As explained in the proof of Theorem \ref{weakilagrangian}, the velocity field ${\mathrm{\mathbf{v}}}$ admits a unique flow. Since $\omega$ is a weak solution, we know also from Theorem \ref{weakilagrangian} that  $\omega(t,X(t,x))=\omega_0(x)$. The specific form of $K_\theta$ then gives us that $\overline\partial \v(t)=\frac{e^{i\theta}}2\,\omega(t)$. Therefore
$$2\|\overline\partial{\v}(t,\cdot)\|_{L^\infty}=  \|\omega(t,\cdot)\|_{L^\infty}=\|\omega_0\|_{L^\infty},$$
since we know from Theorem \ref{weakilagrangian} that the flow $X(t,\cdot)$ preserves Lebesgue null sets. Moreover, from Lemma \ref{LinftyK} (b) we know that $\v(t,z)$ vanishes at $\infty$ like $\frac{C}{|z|}$. Hence,
$$\aligned
\underset{|z|\to\infty}{\limsup}\frac{|\v(t,z)|}{|z|\log\left(e+|z|\right)}&\leq C<+\infty.
\endaligned$$
We are thus in the extended Reimann's situation from \cite{CJMO3}, and so quasiconformality follows, with quasiconformality constant 
$$1\leq \mathcal{K}_t\leq \exp\left(2\int_0^t\|\overline\partial \v(s,\cdot)\|_{L^\infty} ds\right)\leq e^{t\|\omega_0\|_{L^\infty}},$$
this is the claim. 
\end{proof}

 \vspace{.5cm}

Spatial weak derivatives of the flow $\partial X(t,\cdot), \overline\partial X(t,\cdot)$ exist at every time $t$, as a consequence of its quasiconformality. We will now show their absolute continuity in time. To do this, we will use the quasiconformal invariance of the homogeneous Sobolev space $\dot{W}^{1,2}(\C)$.  
 \vspace{.5cm}

\begin{lem}\label{timederivatives}
Let $\omega\in L^\infty([0,T];L^1\cap L^\infty)$ be a weak solution of \eqref{cauchysystem}, and suppose that $\omega_0\in L^\infty$ is compactly supported. Then:
\begin{itemize}
\item[(i)] $D({\mathrm{\mathbf{v}}}(t)\circ X(t))$ belongs to $L^\infty([0,T];L^2 )$.
\item[(ii)] $DX(t,\cdot)-{\mathrm{\mathbf{Id}}}$ belongs to $W^{1,\infty}([0,T];L^2 )$.
\end{itemize}
\end{lem}
\begin{proof}
From Corollary \ref{boundeddispl}, there exists $R>0$ such that $\supp\omega(t)\subset D(0,R)$ whenever $0\leq t\leq T$. Also, since $\v(t)=K_\theta\ast \omega(t)$, we immediately see that $\v(t)$ has distributional derivatives $D\v(t) \in L^2 $, and moreover
$$\aligned
\|D\v(t)\|_{L^2 }^2&=\int_\C(|\partial \v(t)| +|\overline\partial \v(t)|)^2\,dA\\&=\frac14\int_\C (|\cB(\omega(t))|+|\omega(t)|)^2\,dA \leq \int_\C|\omega(t)|^2\,dA\leq \|\omega(t)\|_{L^1}\,\|\omega(t)\|_{L^\infty}.
\endaligned$$
We now use Corollary \ref{qcflows}, from which we know that $|DX(t,\cdot)|^2\leq \cK_t\,\det DX(t,\cdot)$ a.e., and therefore
$$\aligned
\|D (\v(t)\circ &X(t))\|_{L^2 }^2\le \int |D\v(t, X(t,x))|^2 \,|DX(t,x)|^2\,dA(x)\\&\leq \cK_t\int |D\v(t, X(t,x))|^2 \,\det DX(t,x)\,dA(x)=\cK_t\int_\C|D\v(t)|^2\,dA\endaligned$$
thus $(i)$ follows since $\omega\in L^\infty([0,T];L^1\cap L^\infty)$. This allows us to differentiate the integral representation of the flow,
$$X(t,x)=x+\int_0^t \v(s,X(s,x))\,ds$$
which is done by testing against a smooth compactly supported derivative and applying Fubini's Theorem. One immediately gets that
\begin{equation}\label{abscontDX}
DX(t,x)=\Id+\int_0^t D\bigg(\v(s,X(s,x))\bigg)\,ds
\end{equation}
This, combined with $(i)$, shows that $DX(t,\cdot)-\Id$ is the integral in time of a function in $L^\infty([0,T];L^2 )$ whence $DX(t,\cdot)-\Id\in W^{1,\infty}([0,T]; L^2 )$. The claim follows.
\end{proof}

Lemma \ref{timederivatives}, together with the chain rule, legitimates us to take spatial derivatives at \eqref{abscontDX}. One obtains for a.e. $t$ that $\frac{d}{dt} DX(t,x)=D\big(\v(t,X(t,x))\big)$ as $L^2$ functions. In particular,
\begin{equation}\label{chainrule}
\aligned
\frac{d}{dt} \,\partial X(t,z)&= \partial \v(t, X(t,z))\,\partial X(t,z)+\overline\partial \v (t,X(t,z))\,\partial\overline{X}(t,z)
\\
\frac{d}{dt} \,\overline\partial X(t,z)&= \partial \v(t, X(t,z))\,\overline\partial X(t,z)+\overline\partial \v (t,X(t,z))\,\overline{\partial X}(t,z).
\endaligned
\end{equation}
Both equations can be interpreted as the equality of honest $L^2$ functions, and will be useful in the following result, in which we recall a key point in the existence proof in \cite{CS}.  We will write
$$J (t,z)=\det DX (t,z)$$
to denote the Jacobian determinant of $X(t,\cdot)$ at $z$.  

\begin{lem}\label{principal}
Let $\omega\in L^\infty([0,T];L^1\cap L^\infty)$ be a weak solution of \eqref{cauchysystem}, and suppose that $\omega_0\in L^\infty $ is compactly supported. Then  $X(t,\cdot)$ is conformal away from $\supp(\omega_0)$.
\end{lem}

\begin{proof}
After denoting $\mu(t,z)=\frac{\overline\partial X(t,z)}{\partial X(t,z)}$, a bit of manipulation in \eqref{chainrule} gives us that
$$\aligned
\overline\partial \v(t,X (t,z))
=\frac{\partial X(t,z)\,\frac{d}{dt}\overline\partial X(t,z) -\overline\partial X(t,z)\,\frac{d}{dt}\partial X(t,z)}{J(t,z)} 
=\frac{\frac{d}{dt}\mu(t,z)}{1-|\mu(t,z)|^2}\,\frac{\partial X(t,z)}{\overline{\partial X}(t,z)}.
\endaligned$$
On the other hand, from the kernel structure we have
$$
2|\overline\partial \v (t,X (t,z))|= |\omega(t,X (t,z))|=|\omega_0(z)|.
$$
Thus
$$
\aligned
\frac{ \frac{d}{dt}|\mu (t,z)|}{1-|\mu (t,z)|^2}
&\leq\frac{\left|\frac{d}{dt}\mu (t,z)\right|}{1-|\mu (t,z)|^2}=\frac12\,|\omega_0(z)|
\endaligned$$
Now, given any time $t>0$, we can integrate on $(0,t)$ the above inequality to obtain that
\begin{equation}\label{pointw}
\log\left(\frac{1+|\mu(t,z)|}{1-|\mu(t,z)|}\right)\leq t\,|\omega_0(z)|,
\end{equation}
since $X(0,z)=z$ implies $\mu(0,z)=0$. Now, if $\omega_0(z)=0$ then necessarily $\mu(t,z)=0$ and hence $\overline\partial X(t,z)=0$, so the claim follows.  
\end{proof}

\section{The energy inequality}\label{energy}

We consider now a compactly supported initial datum $\omega_0\in L^\infty$, and suppose that $\omega_1,\omega_2\in L^\infty([0,T]; L^1\cap L^\infty)$ are weak solutions of \eqref{cauchysystem} with initial datum $\omega_0$. From Theorem \ref{weakilagrangian}, we know each $\omega_i$ is also a Lagrangian solution, and thus each $\omega_i$ comes together with a velocity field $\v_i=K_\theta\ast\omega_i$ whose flow $X_i$ satisfies that $\omega_i(t,X_i(t,x))=\omega_0(x)$ at a.e. $x\in\C$.  Note that for each $t\in[0,T]$ both flows $X_i(t)$ are $\cK_t$-quasiconformal, and in both cases $\cK_t\leq e^{t\|\omega_0\|_{L^\infty}}$ by Corollary~\ref{qcflows}. Being on the plane, quasiconformal mappings are better understood via their Beltrami equation, so we set $\mu_i(t)=\frac{\overline\partial X_i(t)}{\partial X_i(t)} $, which is supported in $\supp(\omega_0)$ by Lemma~\ref{principal}. By the quasiconformality of $X_i$ for every $t$, one has $ \|\mu_i(t)\|_{L^\infty}<1$. We will need uniform ellipticity in the forthcoming  arguments, and thus we fix a constant $k_0<1$ and declare $\tau = \frac1{\|\omega_0\|_{L^\infty}}\,\log\left(\frac{1+k_0}{1-k_0}\right).$
Thus it holds that
\begin{equation}\label{elliptbound}
\frac{e^{t\|\omega_0\|_{L^\infty}}-1}{e^{t\|\omega_0\|_{L^\infty}}+1}\leq k_0\qquad\text{whenever }0\leq t\leq \min\{T,\tau\}.
\end{equation}
By Corollary~\ref{qcflows}, this gives that $\|\mu_i(t)\|_{\infty}\le k_0$ for each $t\in[0,\tau]$. Together with $k_0$, we also choose an exponent $p_0$ such that
\begin{equation}\label{ellipticp0}
2<p_0<1+\frac1{k_0}.
\end{equation}
Within that range, the Beltrami operator $\Id-\mu_i(t)\cB$ is boundedly invertible in $L^{p_0}$ (see \cite{AIS} or also \cite{AFBeurling}). Now that the parameters $\tau,k_0,p_0$ are fixed we carry on with our search of the right ODE for the energy. It will be also important to fix an auxiliary exponent $\gamma_0$ defined by 
\begin{equation}\label{gamma0def}
\frac2{p_0}+\frac2{\gamma_0}=1.
\end{equation}
The following two lemmas will be used in deducing the energy inequality. 
\begin{lem}\label{uniformqcbound}
Let $\omega\in L^\infty([0,T];L^1\cap L^\infty)$ be a weak solution of $\eqref{cauchysystem}$. We have
$$\|\overline\partial X(t)\|_{L^{p_0}}+\|\partial X(t)-1\|_{L^{p_0}}\leq \frac{|D|^{\frac1{p_0}}\,p_0}{1+\frac1{k_0}-p_0}$$
for any disk $D\supset \supp(\omega_0)$ and $0\leq t\leq\min\{T,\tau\}$.
 \end{lem}
\begin{proof}
We know from Corollary \ref{qcflows} and Lemma \ref{principal} that $X(t)$ is $\cK_t$-quasiconformal on $\C$ and conformal away from $\supp(\omega_0)$. Thus, for any fixed $t\in[0,T]$, if we set $\mu(t)=\frac{\overline\partial X(t)}{\partial X(t)}$, then $\mu(t)$ is compactly supported, essentially bounded with $\|\mu(t)\|_{L^\infty}\leq\frac{\cK_t-1}{\cK_t+1}$, and moreover the equation 
\begin{equation}\label{beltrami}
\overline\partial X(t)=\mu(t)\,\partial X(t)=\mu(t)\,\partial (X(t)-\Id)+\mu(t)
\end{equation}
holds almost everywhere. From Corollary \ref{boundeddispl} (c), we can infer also that $X(t,z)-z=\sum_{n\geq 1}\frac{a_n(t)}{z^n}$ for large values of $z$, thus granting $\partial (X(t)-\Id)\in L^p$ for some $p>2$. Thus \eqref{beltrami} gives also that $\overline\partial (X(t)-\Id)\in L^p$. This allows us to infer that $\partial X(t)-1=\cB(\overline\partial X(t))$, whence \eqref{beltrami} reads as
$$(\Id-\mu(t)\,\cB)(\overline\partial X(t))=\mu(t).$$
Condition \eqref{elliptbound} and our choice of $p_0$ guarantee that the operator $\mu(t)\cB$ is bounded in $L^{p_0}$ with norm
$$\|\mu(t)\cB\|_{L^{p_0}\to L^{p_0}}\leq \frac{e^{t\|\omega_0\|_{L^\infty}}-1}{e^{t\|\omega_0\|_{L^\infty}}+1}\,(p_0-1)\leq k_0(p_0-1)<1$$ 
at least for values $0\leq t\leq \min\{T,\tau\}$. Thus, a Neumann series argument proves the invertibility of $\Id-\mu(t)\cB$ on $L^{p_0}$ for each $t\in[0,\tau]$, and moreover
$$\|\overline\partial X(t)\|_{L^{p_0}}\leq \frac{\|\mu(t)\|_{L^{p_0}}}{1-k_0(p_0-1)} \leq \frac{k_0|D|^\frac1{p_0}}{1-k_0(p_0-1)}.$$
Next, since $\partial X(t)-1=\cB(\overline\partial X(t))$, the boundedness of $\cB$ on $L^{p_0}$ implies that
$$\| \partial X(t)-1\|_{L^{p_0}}\leq  \frac{k_0|D|^{\frac1{p_0}}\,(p_0-1)}{1-k_0(p_0-1)}.$$
This finishes the proof.
\end{proof}
\begin{lem}\label{linearinq}
There is a constant $C=C(T, k_0, p_0,\|\omega_0\|_{L^\infty} )$ such that
$$\|\partial{\mathrm{\mathbf{v}}}(t, X(t))\|_{L^q(D)}\leq C\,q|D|^{1/q} ,$$
for each $q>\gamma_0$ and each $t\in[0,\min\{T,\tau\}]$, and for every disk $D\supset\supp(\omega_0)$.  
\end{lem}

\begin{proof}
By H\"older's inequality,
$$\aligned
\|\partial{\mathrm{\mathbf{v}}}(t,X(t))\|_{L^q(D)}^q
&=\int_{X(t,D)}  |\partial\v(t,x)|^q\,\det(D(X^{-1}(t,x)))\,dA(x)\\
&\leq \left(\int_{X(t,D)}  |\partial\v(t,x)|^\frac{\gamma_0q}2\,dA(x)\right)^\frac{2}{\gamma_0}\,\left(\int_{X(t,D)}\det(DX^{-1}(t,x))^\frac{p_0}2\,dA(x)\right)^{\frac{2}{p_0}}.
\endaligned$$
The last integral above is easily handled with the help of the reverse H\"older inequality from \cite[Theorem 13.2.7]{AIM}, which can be applied since from \eqref{ellipticp0} we have $1-\frac{p_0}2>\frac12(1-\frac1{k_0})$. Indeed, 
$$
\aligned
\int_{X(t,D)}\det ( D X^{-1}(t,\cdot))^\frac{p_0}2
&= \int_{D}\det ( D X(t,\cdot))^{1-\frac{p_0}2} \\
&\leq C(\cK_t)\,|D| \,\left(\fint_{D}\det ( D X(t,\cdot))\right)^{1-\frac{p_0}{2}}= C(\cK_t)\,|D|^{\frac{p_0}2}\,|X(t,D)|^\frac{2-p_0}{2} 
\endaligned
$$
with a constant $C(\cK_t)\geq 1$ that increases with $\cK_t$. Concerning the first term, and using that $\frac{\gamma_0q}{2}>1$, 
$$
\aligned
\left(\int_{X(t,D)}  |\partial\v(t,x)|^\frac{\gamma_0q}2\,dA(x)\right)^\frac{2}{\gamma_0q}
&=\frac12\left(\int_{X(t,D)}  |\cB(\omega(t))|^\frac{\gamma_0q}2\,dA(x)\right)^\frac{2}{\gamma_0q}\\
&\leq \frac{\frac{\gamma_0q}2-1}{2}\left(\int_{X(t,D)}  | \omega(t)|^\frac{\gamma_0q}2\,dA(x)\right)^\frac{2}{\gamma_0q}\\
&\leq \frac{\frac{\gamma_0q}2-1}{2}\,\|\omega_0\|_{L^\infty}\,|X(t,D)|^\frac{2}{\gamma_0q}. 
\endaligned$$
Putting things together, and using \eqref{gamma0def}, we get the desired estimate.
\end{proof}



\noindent
We now introduce
$E(t)=\|\overline\partial X_1(t)-\overline\partial X_2(t)\|_{L^2}.$
Lemma \ref{timederivatives} asserts that $E(t)$ is well-defined, and moreover $\overline\partial X_i \in W^{1,1}(0,T;L^2(\C))$ so that $E$ is absolutely continuous with derivative
$$\frac{d}{dt}\,E(t)\leq\left\|\frac{d}{dt}\,\bigg(\overline\partial X_1(t)- \overline\partial X_2(t)\bigg)\right\|_{L^2}$$
for almost every $t\in[0,T]$.  

\begin{theo}[The energy inequality]\label{thm:energyinequality}
Let $\tau,\gamma_0$  be defined as above and $q>\gamma_0$. Then there are constants $C_1,C_2$depending on $(k_0,\omega_0)$ such that
\begin{equation}\label{eq:Epower}
 E'(t)\leq C_1\,q\,E(t)^{1-\frac{\gamma_0}{q}}+ C_2\,E(t),\hspace{1cm}a.e. \,\,0\leq t\leq \min\{T,\tau\}.
\end{equation}
\end{theo}

\begin{proof}

We now call $H(t)=X_1(t)-X_2(t)$. Lemmas \ref{timederivatives} and \ref{principal} allow us to claim that $H(t)\in \dot{W}^{1,2}$ for each $t$ and moreover $DH \in W^{1,\infty}([0,T]; L^2 )$. In particular, we have the following identities
$$\aligned
\partial H(t)&=\cB(\overline\partial H(t))\\
 \dot{\partial H}(t) &=\cB(\dot{\overline\partial H}(t)) 
\endaligned$$
of $L^2$ maps, where the dot denotes time derivative. We now use them to get an ODE for $\overline\partial H(t)$, from which we will obtain an ODE for $E(t)$. We first subtract the equations for $\dot{\overline\partial X_1}(t)$ and $\dot{\overline\partial X_2}(t)$ from \eqref{chainrule}.  Inserting the definition of $\mu_1$ in the evolution equations for $X_i(t)$
we obtain that 
$$\aligned
\dot{\overline\partial H}(t)-\mu_1(t)\dot{ \partial H}(t)
&=\partial\v_1(t,X_1(t))\,\left(\overline\partial X_1(t)-\mu_1(t)\partial X_1(t)\right)-\partial\v_2(t,X_2(t))\left(\overline\partial X_2(t)-\mu_1(t)\partial X_2(t)\right)\\&+\frac{\omega_0}{2}\,\left(\overline{\partial H}(t)-\mu_1(t)\overline{\overline\partial H}(t)\right)
\endaligned$$
The first term on the right hand side vanishes since $\overline\partial X_1(t)-\mu_1(t)\partial X_1(t)=0$ a.e. Concerning the second term, we note that
$$
-\left(\overline\partial X_2(t)-\mu_1(t)\partial X_2(t)\right)=\overline\partial H(t)-\mu_1(t)\partial H(t)$$
while for the third we simply observe that
$$\overline{\partial H}(t)-\mu_1(t)\overline{\overline\partial H}(t)=\overline{\partial H}(t)-\mu_1(t) \partial \overline{H}(t).$$
Now we rewrite the Beltrami-type equation in terms of the Beurling transform. At almost every time, it holds that
\begin{equation}\label{fundamental}\aligned
(\Id-\mu_1(t)\cB)(
\dot{\overline\partial H}(t)) 
&= \partial\v_2(t,X_2(t))\left(\Id-\mu_1(t)\cB\right)(\overline\partial H(t)) \\&+\frac{e^{i\theta}\,\omega_0}{2}\,\left(\Id-\mu_1(t)\cB\right)(\overline{\partial H}(t)).
\endaligned
\end{equation}
By Lemma \ref{principal}, all members in the above equation have uniformly bounded support, say, on a disk $D\supset \supp(\omega_0)$, and therefore for each such term the $L^p$ norm can be taken only on $D$, that is $L^p=L^p(D)$.  We first use the left hand side to control $E'(t)$ since we can use the Neumann series to invert the Beltrami Operator. Namely,  
$(I-\mu_1 \cB)^{-1}=\sum (\mu_1 \cB)^n$,
which yields the bound
\[
\| (\Id-\mu_1(t)\cB)^{-1}\|_{L^2\to L^2}\leq \frac1{1-k_0}.\]
for $0\le \min\{\tau,T\}$.
This allows us to bound from below the term on the left at \eqref{fundamental} as follows,
$$E'(t)\leq\|\dot{\overline\partial H}(t)\|_{L^2}\leq\frac1{1-k_0}\,\|(\Id-\mu_1(t)\cB)(
\dot{\overline\partial H}(t)) \|_{L^2},$$
for all times $0\leq t\leq \min\{T,\tau\}$. Concerning the second term on the right hand side of \eqref{fundamental}, we have for every $t\geq 0$ the direct estimate
$$
\left\|\frac{\omega_0}{2}\,\left(\Id-\mu_1(t)\cB\right)(\overline{\partial H}(t))\right\|_{L^2}\leq\\|\omega_0\|_{L^\infty}\,E(t).
$$
Finally, in the first term on the right hand side of \eqref{fundamental} we take $q,\tilde{q}>2$ such that $\frac12=\frac1q+\frac1{\tilde{q}}$ and use H\"older's inequality on $D$, and obtain
\begin{equation}\label{last2terms}
\aligned
\|\partial\v_2(t,X_2(t))\left(\Id-\mu_1(t)\cB\right)(\overline\partial H(t))\|_{L^2(D)}\leq \|\partial\v_2(t,X_2(t))\|_{L^q(D)}\,\|\left(\Id-\mu_1(t)\cB\right)(\overline\partial H(t))\|_{L^{\tilde{q}}(D)}
\endaligned
\end{equation}
For the second term on the right hand side of \eqref{last2terms}, we observe that
\begin{equation}\label{firstof2terms}
\|\left(\Id-\mu_1(t)\cB\right)(\overline\partial H(t))\|_{L^{\tilde{q}}(D)}\leq (1+k_{0}(\tilde{q}-1))\|\overline\partial H(t)\|_{L^{\tilde{q}}}
\end{equation}
and now we bound the last factor on the right hand side above by combining interpolating estimates with the result below. We restrict our attention  to values of $q$ for which $2<\tilde{q}<p_0$. Namely, by writing
$$\frac1{\tilde{q}}=\frac{1-\theta}{2}+\frac{ \theta}{p_0} $$
one gets necessarily that $\left(\frac12-\frac1{p_0}\right)\theta=\frac1q$. Since $\frac1{p_0}+\frac1{\gamma_0}=\frac12$, this actually means that $\theta=\frac{\gamma_0}{q}$, and therefore the only requirement to have $2<\tilde{q}<p_0$ is to choose $q>\gamma_0$. If so, the interpolating estimate we get is
\begin{equation}\label{interpolating}
\|\overline\partial H\|_{L^{\tilde{q}}}\leq \|\overline\partial H\|_{L^2}^{1-\theta}\,\|\overline\partial H\|_{L^{p_0}}^{\theta}
\le E(t)^{1-\frac{\gamma_0}{q}}\,\left(\frac{2|D|^\frac1{p_0}\,p_0}{1+\frac1{k_0}-p_0}\right)^{\gamma_0/q}
\end{equation}
after using Lemma \ref{uniformqcbound}. 
We now plug the estimate from Lemma \ref{linearinq} into \eqref{last2terms}, and combine it with the interpolating estimates \eqref{interpolating}. We obtain for any $q>\gamma_0$ that
$$ 
\|\partial\v_2(t,X_2(t))\left(\Id-\mu_1(t)\cB\right)(\overline\partial H(t))\|_{L^2(D)}\leq 
Cq|D|^{1+\frac1{p_0}}\,(1+k_{0}(\tilde{q}-1))\,E(t)^{1-\frac{\gamma_0}{q}}$$
with a constant $C$ that may depend on $p_0,\|\omega_0\|_{L^\infty}$ but not on $D$ or $q$.  
Hence \eqref{fundamental} becomes \eqref{eq:Epower} and the energy inequality is proved.
\end{proof}

\section{ Proof of Theorem 1}\label{proofofuniqueness}

After having proven \eqref{eq:Epower}, the proof of uniqueness follows the standard Gronwall-type argument, which we specify below for the reader's convenience. Choose a fixed constant $A$ satisfying 
\begin{equation}
A\ge\max\left\{1,\ e^2\sup_{0\le t\le T_*}E(t)\right\},
\label{defA}
\end{equation}
where we wrote $T_\ast=\min\{T,\tau\}$, and set $\eta(t)=\frac{E(t)}{A}.$ By construction 
$$
0\le \eta(t)\le e^{-2},\qquad 0\le t\le T_*.
$$
Now,  dividing \eqref{eq:Epower} by $A$ we get the differential inequality 
$$
\eta'(t)\le  C_1\,q\,\eta(t)^{1-\frac{\gamma_0}{q}}+C_2\,\eta(t), \text{ for almost every } t\in[0,T_{*}] \text{ and } q>\gamma_0.
$$
After excluding the null set of times where $E$ is not differentiable, we have a common set of full measure on which all the inequalities hold simultaneously. So we now fix a time $t$ in this common set with $\eta(t)>0$ and define
$$ 
\ell(t)=\log\frac{1}{\eta(t)}. 
$$
Then $\ell(t)\geq 2$ and moreover
$$
q\,\eta(t)^{1-\frac{\gamma_0}{q}}=\eta(t)\,q\,\eta(t)^{-\frac{\gamma_0}{q}}
=\eta(t)\,q\exp\left(\frac{\gamma_0\,\ell(t)}{q}\right).
$$
We now choose the integer $q= [\gamma_0\,\ell(t)]+1,$ and obtain therefore that $q\ge 2\gamma_0$ so in particular $q>\gamma_0$. Moreover, $\exp\left(\frac{\gamma_0\,\ell(t)}{q}\right)\le e$ and
$$q\le \gamma_0\,\ell(t)+1\le\left(\gamma_0+\frac{1}{2}\right)\ell(t) $$
Consequently, 
$$
q\eta(t)^{1-\frac{\gamma_0}{q}}\le e\,\ell(t)\,\eta(t)\,\left(\gamma_0+\frac{1}{2}\right).
$$
Since $\ell(t)\ge 2,$ the remaining term $C_2\eta(t)$ is bounded by a constant times $\eta(t) \,\ell(t)$ and so enlarging the constants if needed we end up getting that
$$
\eta'(t)\le C\,\eta(t)\,\log\frac1{\eta(t)}
$$
which after multiplying by \(A\) gives
\begin{equation}\label{eq:Osgood}
 E'(t)\le C\,E(t)\,\log\left(\frac{A}{E(t)}\right),\qquad a.e.t\in[0,T_\ast] .
\end{equation}
We now prove that \(E\) cannot be positive. To this end, we define for $r>0$ the function 
$$
\varphi_A(r)=r\log\left(\frac{A}{r}\right) ,\qquad \varphi_A(0)=0,
$$
which is continuous on $[0,\infty)$. For \(0<r\le A/e\) it is also differentiable and
\[
\varphi_A'(r)=\log\left(\frac Ar\right)-1\ge0,
\]
so \(\varphi_A\) is increasing on \([0,A/e]\). 
For the next arguments we will have $E$ in the denominator but it could be zero and for this reason we add a positive constant before dividing. Take $0<\delta<A/e-A/e^2$ and define $E_\delta(t)=E(t)+\delta.$ Since $0 \le E(t)\le A/e^2,$ we have $ 0<E_\delta(t)<A/e$ on the whole interval. The function \(E_\delta\) is absolutely continuous, and the monotonicity of \(\varphi_A\) gives
$$
E_\delta'=E'\le C\varphi_A(E)\le C\varphi_A(E_\delta).
$$
Define next
$$
L_\delta(t)=\log\left(\frac{A}{E_\delta(t)}\right).
$$
For fixed \(\delta>0\), the denominator stays at least \(\delta\), so the chain rule for the composition of a Lipschitz function with an absolutely continuous function applies and we get 
$$
L_\delta'(t)=-\frac{E_\delta'(t)}{E_\delta(t)}\ge -C \frac{\varphi_{A}(E_{\delta})}{E_{\delta}}
\ge-CL_\delta(t).
$$
Moreover we obtain almost everywhere the inequality
$$
\frac{d}{dt}\bigl(e^{Ct}L_\delta(t)\bigr)
=e^{Ct}\bigl(L_\delta'(t)+CL_\delta(t)\bigr)\ge0
$$
By absolute continuity we can do integration between any \(0\le t_1\le t_2\le T_*\), obtaining
$$
L_\delta(t_2)\ge e^{-C(t_2-t_1)}L_\delta(t_1).
$$
After exponentiating and rearranging, this gives
$$
E(t_2)+\delta
\le A\left(\frac{E(t_1)+\delta}{A}\right)^{e^{-C(t_2-t_1)}}.
$$
The next argument is to remove the regularization given by $\delta$.  Since the constants $A$ and $C$ do not depend on $\delta$, and the power on the right hand side is positive, we can let $\delta\searrow 0$ to obtain  the inequality 
\begin{equation}\label{eq:integrated}
E(t_2)\le A\left(\frac{E(t_1)}A\right)^{\exp(-C(t_2-t_1))},
\qquad 0\le t_1\le t_2\le T_*.
\end{equation}
Consequently, if $E(t_1)=0$ we have $E(t_2)=0$ for every $t_2\in[0,T_{*}].$

\vspace{1cm}
To conclude the proof, you will argue as follows. Initially $X_1(0)=X_2(0)=\Id,$ so that $E(0)=0$. Applying inequality \eqref{eq:integrated} with $t_1=0$ and $t_2=t$ gives $E(t)=0$ for every $t\in[0,T_{*}]$. Consequently,
$$
\overline\partial(X_1-X_2)(t)=0\quad\text{in }L^2.
$$
and therefore
$$
\partial(X_1-X_2)(t)=\cB\bigl(\overline\partial(X_1-X_2)(t)\bigr)=0\quad\text{in }L^2. 
$$
Thus all spatial derivatives of $X_1(t)-X_2(t)$ vanish almost everywhere, and thus $X_1(t)-X_2(t)$ must be a constant. Corollary \ref{boundeddispl} (c) now tells us that the only viable option is $X_1(t)=X_2(t)$. Using now Theorem  \ref{weakilagrangian}, we obtain 
$$
\omega_1(t)=\omega_0\circ X_1(t)^{-1}
=\omega_0\circ X_2(t)^{-1}=\omega_2(t)
$$
throughout the interval $[0,T_{*}]$.

\vspace{.5cm}
Finally, let \(t_0\) be the endpoint of an interval on which equality has been established. Next, use the common datum \(\omega_1(t_0)=\omega_2(t_0)\), still bounded and compactly supported, and moreover $\|\omega_1(t_0)\|_{L^\infty}=\|\omega_2(t_0)\|_{L^\infty}=\|\omega_0\|_{L^\infty}$ to restart the problem. Then one gets two flows \(X_i(t,t_0,\cdot)\) which are equal to the identity at time $t_0$, and so their derivative comparison again starts from zero.

The distortion estimate now involves the elapsed time \(t-t_0\). Thus the same threshold \(\tau\), which depends only on $\|\omega_0\|_{L^\infty}$ and $k_0$, gives the required small Beltrami coefficients on \([t_0,\min\{t_0+\tau,T\}]\). Repeating the above arguments with this common datum proves equality on the next interval. The constants $A$ and $C$ may change with the support and velocity bounds, but they remain finite, which is all that the zero-energy conclusion requires. Finitely many consecutive intervals cover $[0,T]$ and it completes the proof of Theorem~\ref{main1}.

\section{Proof of Theorem \ref{main2}}\label{stability}

We consider now a pair $\omega_{01},\omega_{02}\in L^\infty$, both compactly supported in $D=D(0,R_0)$. For $i=1,2$, let $\omega_i\in L^\infty([0,T]; L^1\cap L^\infty)$ be the unique weak solution of \eqref{cauchysystem} with initial datum $\omega_{0i}$, which exists and is unique due to Theorem \ref{main1}. From Theorem \ref{weakilagrangian}, we know each $\omega_i$ is also a Lagrangian solution, and thus each $\omega_i$ comes together with a velocity field $\v_i=K_\theta\ast\omega_i$ whose flow $X_i$ satisfies that $\omega_i(t,X_i(t,x))=\omega_{0i}(x)$ at a.e. $x\in\C$.  Calling $M=\max\{\|\omega_{01}\|_{L^\infty},\|\omega_{02}\|_{L^\infty}\}$, for each $t\in[0,T]$ both flows $X_i(t)$ are $\cK_t$-quasiconformal with  $\cK_t\leq e^{tM}$. Similarly to what we did in Section \ref{energy}, we introduce a new ellipticity bound $k_0\in(0,1)$ and choose a time $\tau$ such that $\cK_t\leq \frac{1+k_0}{1-k_0}$ for each $t\in[0,\min\{T,\tau\}]$, namely
$$\tau=\frac1{M}\,\log\left(\frac{1+k_0}{1-k_0}\right).$$
Next, we call again
$$E(t)=\|\overline\partial X_1(t)-\overline\partial X_2(t)\|_{L^2}.$$
As before, Lemma \ref{timederivatives} guarantees we have enough regularity to take time derivatives.  Calling $H(t)=X_1(t)-X_2(t)$, and $\mu_1(t)=\frac{\overline\partial X_1(t)}{\partial X_1(t)}$, we now have instead of \eqref{fundamental} the following identity
\begin{equation}\label{fundamental2}\aligned
(\Id-\mu_1(t)\cB)(
\dot{\overline\partial H}(t)) 
&= \partial\v_2(t,X_2(t))\left(\Id-\mu_1(t)\cB\right)(\overline\partial H(t)) \\&+\frac{e^{i\theta}(\omega_{01}-\omega_{02})}{2}\,\left(\Id-\mu_1(t)\cB\right)(\overline{\partial X_1}(t))+\frac{e^{i\theta}\,\omega_{02}}{2}\,\left(\Id-\mu_1(t)\cB\right)(\overline{\partial H}(t)),
\endaligned
\end{equation}
for a.e. $t\in[0,T]$. After taking $L^2(D)$ norms, the invertibility bounds for $\Id-\mu_1(t)\cB$ force us to look only at $t\leq \tau$. Having this in mind, the term on the left hand side and the first and last terms on the right hand side are bounded exactly as we did in Section \ref{energy}. Concerning the central term on the right hand side, we observe that
$$\aligned
\left\|\frac{\omega_{01}-\omega_{02}}{2}\,\left(\Id-\mu_1(t)\cB\right)(\overline{\partial X_1}(t))\right\|_{L^2(D)}
&\leq \frac12\,\|\omega_{01}-\omega_{02}\|_{L^{\gamma_0}}\,\|\left(\Id-\mu_1(t)\cB\right)(\overline{\partial X_1}(t))\|_{L^{p_0}(D)}\\
&\leq C(k_0,p_0,D)\,\|\omega_{01}-\omega_{02}\|_{L^{\gamma_0}} 
\endaligned$$
whence \eqref{fundamental2} becomes
\begin{equation}\label{Epower2}
 E'(t)\leq C_1\,q\,E(t)^{1-\frac{\gamma_0}{q}}+ C_2\,E(t)+C_3\,\|\omega_{01}-\omega_{02}\|_{L^{\gamma_0}} \hspace{1cm}a.e. \,\,0\leq t\leq \min\{T,\tau\}.
\end{equation}
for any $q>\gamma_0$, with $C_1,C_2,C_3$ depending on $M, k_0,p_0$ and $|D|$, but not on $q$. By a similar argument, this equation reduces to
\begin{equation}\label{simpeqforE}
E'(t)\leq C\, \left(E(t)\,\log\left(\frac{A}{E(t)}\right)+ \|\omega_{01}-\omega_{02}\|_{L^{\gamma_0}}\right), \hspace{1cm}a.e. \,\,0\leq t\leq \min\{T,\tau\}
\end{equation}
for some conveniently chosen constant $A\geq 1$.  Calling then $\delta_0=\|\omega_{01}-\omega_{02}\|_{L^{\gamma_0}}$ and $Z(t) := E(t) + \delta_0$, one gets $ E'(t)=Z'(t) \le C\left(E\log(A/E)+\delta_0\right).$ Choosing a larger $A$ if needed, one can assume that $E(t), Z(t) \in [0,\frac{A}e]$. It is easy to verify that the function $f(x)=x\log(A/x)$ is increasing on this interval and hence
\begin{equation}
E\log(A/E)\le Z\log(A/Z).
\label{firstinE}
\end{equation}
Also we have
\begin{equation}
\delta_0 \le Z \le Z\log(A/Z).
\label{secondinE}
\end{equation}
Inserting \eqref{firstinE} and \eqref{secondinE} into \eqref{simpeqforE} we have the following differential inequality for $Z$,
$$\begin{cases}Z' \le C Z \log(A/Z),\\
Z(0)= \delta_0.
\end{cases}$$
We can solve it and obtain
$$E(t) \le Z(t) \le A^{1-\theta(t)}\delta_0^{\theta(t)},$$
where $\theta(t):=\mathrm{e}^{-Ct}.$ Theorem \ref{main2} follows by interpolation between $L^1$ and $L^{\infty}$.

\end{document}